# A Geometric-Probabilistic problem about the lengths of the segments intersected in straights that randomly cut a triangle.

*Jesús Álvarez Lobo.*
*Spain.*

**Abstract.** If a line cuts randomly two sides of a triangle, the length of the segment determined by the points of intersection is also random. The object of this study, applied to a particular case, is to calculate the probability that the length of such segment is greater than a certain value.

*Let ABC be an isosceles triangle, with $\overline{AB} = \overline{CB}$ and $\overline{OB} = \overline{AC}$, being O the midpoint of $\overline{AC}$ (ie, $\overline{OB}$ is the height relative to the side $\overline{AC}$).*

*Through a randomly chosen point P on $\overline{AC}$ is drawn a straight r with also randomly chosen slope. Let Q and R be the points where r intersects $\overline{AB}$ and $\overline{CB}$, respectively.*

*Calculate the probability for the following inequalities:*

$$\boxed{\overline{PQ} > \overline{AC} \ or \ \overline{PR} > \overline{AC}} \tag{1}$$

Let us draw an arc of radius $\overline{AC}$ with center $P$. Let $P$ and $Q$ be the intersection points of this arc with the sides $\overline{AB}$ y $\overline{CB}$, respectively, as shown in the following picture, with the triangle represented in an orthonormal coordinate system, with origin at $O$, x-axis (abscissas) in the direction $OA$ and y-axis (ordinates) in the direction $OB$.

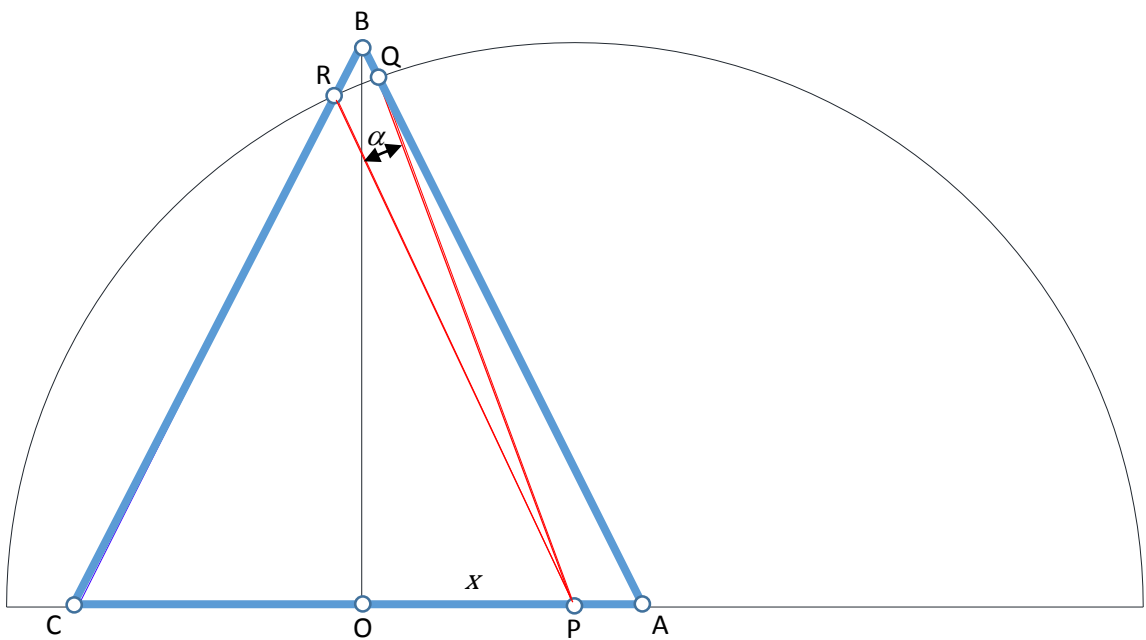

Clearly, all the straight lines of the bundle with vertex P in $\overline{AC}$ intersect the sides $\overline{AB}$ or $\overline{CB}$, and all the lines of the sub-bundle inner to the angle $\alpha = \widehat{QPR}$, and only them, satisfies (1).

Since $x$ and $\alpha$ are **continuous** random variables **uniformly distributed**, for a differential of length $dx$ in $\overline{AC}$, the probability that the condition (1) is satisfied will be

$$dp = \frac{\alpha}{\pi} dx \tag{2}$$

and therefore, the probability that the inequalities (1) are satisfied for a randomly chosen point in $\overline{AC}$ will be

$$p = \frac{1}{\pi} \int_{-\frac{1}{2}}^{\frac{1}{2}} \alpha(x)\, dx \tag{3}$$

where $\alpha(x)$ is the function relating the angle $\alpha$ with the abscissa $x$.

We have used the following facts:

- In (2):
  - In an infinitesimal length, $dx$, the limit angle $\alpha$ is **constant**.
  - The slope of the secant line is **independent** of the abscissa $x$.

- In (3):
  - The required probability $p$ is obtained by Riemann integration of the **probability density function** $\alpha(x)$ in the symmetric interval $\overline{AC} = \left[-\frac{1}{2}, \frac{1}{2}\right]$.

The limit angle $\alpha$ can be expressed in radians as:

$$\alpha = \pi - \widehat{QPA} - \widehat{CPR} \tag{4}$$

But,

$$\overset{\triangle}{QPA}: \quad \widehat{QPA} = \pi - \hat{A} - \hat{Q} \tag{5}$$

$$\overset{\triangle}{CPR}: \quad \widehat{CPR} = \pi - \hat{C} - \hat{R} \tag{6}$$

So,

$$\alpha = \hat{Q} + \hat{R} + \hat{A} + \hat{C} - \pi \tag{7}$$

Perhaps the easiest way to define $\alpha$ as a function of $x$ is trigonometrically:

$$\overset{\triangle}{QPA}: \quad \frac{\sin\hat{A}}{\overline{PQ}} = \frac{\sin\hat{Q}}{\overline{AP}} \Rightarrow \sin\hat{Q} = \frac{\overline{AP}}{\overline{PQ}} \sin\hat{A} \tag{8}$$

$$\overset{\triangle}{CPR}: \quad \frac{\sin\hat{C}}{\overline{PR}} = \frac{\sin\hat{R}}{\overline{CP}} \Rightarrow \sin\hat{R} = \frac{\overline{CP}}{\overline{PR}} \sin\hat{C} \tag{9}$$

But $\hat{C} = \hat{A}$ and $tan\hat{A} = 2$. Moreover, without loss of generality, we can assume that

$$\overline{AC} = \overline{OB} = 1 \tag{10}$$

Hence,

$$sin\hat{A} = sin\hat{C} = \frac{2}{\sqrt{5}} \tag{11}$$

Applying (10) and (11) to (8) and (9), these reduce to

$$sin\hat{Q} = \frac{1-2x}{\sqrt{5}} \tag{8'}$$

$$sin\hat{R} = \frac{1+2x}{\sqrt{5}} \tag{9'}$$

Substituting in (7) this results we get the **probability density function** for the random variable $\alpha$:

$$\alpha(x) = arcsin\left(\frac{1-2x}{\sqrt{5}}\right) + arcsin\left(\frac{1+2x}{\sqrt{5}}\right) + 2arctan(2) - \pi \tag{12}$$

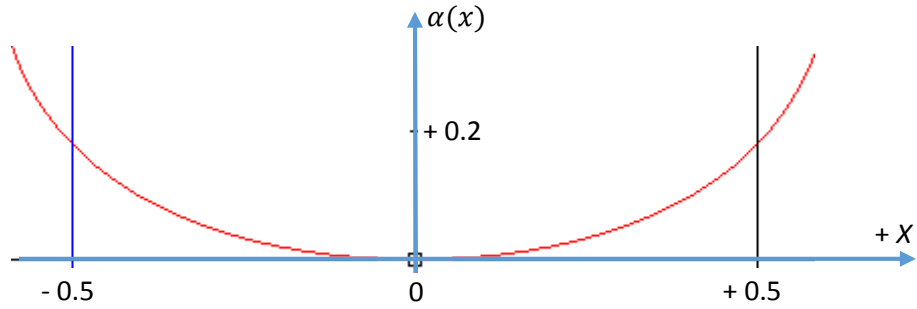

Now, substituting in (3) the result given in (12), we obtain,

$$p = \frac{1}{\pi}\left[\int_{-\frac{1}{2}}^{\frac{1}{2}} arcsin\left(\frac{1-2x}{\sqrt{5}}\right)dx + \int_{-\frac{1}{2}}^{\frac{1}{2}} arcsin\left(\frac{1+2x}{\sqrt{5}}\right)dx\right] + \frac{2}{\pi}arctan(2) - 1 \tag{13}$$

These integrals (in indefinite form) can be solved by *integration by parts*. Let

$$I_1 = \int arcsin\left(\frac{1-2x}{\sqrt{5}}\right) dx \tag{14}$$

$$I_2 = \int arcsin\left(\frac{1+2x}{\sqrt{5}}\right) dx \tag{15}$$

$$I_1 = \int u\,dv \begin{cases} u = arcsin\left(\frac{1-2x}{\sqrt{5}}\right) \\ dv = dx \end{cases} \Rightarrow \begin{cases} du = \dfrac{-\frac{2}{\sqrt{5}}}{\sqrt{1-\left(\frac{1-2x}{\sqrt{5}}\right)^2}} dx \\ v = x \end{cases}$$

And applying the **formula of integration by parts**,

$$I_1 = uv - \int v\,du = x\,arcsin\left(\frac{1-2x}{\sqrt{5}}\right) - \int \frac{-\frac{2}{\sqrt{5}}x}{\sqrt{1-\left(\frac{1-2x}{\sqrt{5}}\right)^2}}\,dx \tag{16}$$

Let

$$I_3 = \int \frac{-\frac{2}{\sqrt{5}}x}{\sqrt{1-\left(\frac{1-2x}{\sqrt{5}}\right)^2}}\,dx \tag{17}$$

After simplifying the sub-integral expression, through the elementary transformations shown below, $I_3$ is reduced to two *quasi-immediate integrals* (reducible to immediate integrals by simple adjustment of constants). Omitting integration constants, for simplicity:

$$I_3 = \int \frac{-x\,dx}{\sqrt{-x^2+x+1}} = \int \frac{-2x+1-1}{2\sqrt{-x^2+x+1}}\,dx = \int \frac{-2x+1}{2\sqrt{-x^2+x+1}}\,dx - \int \frac{dx}{2\sqrt{-x^2+x+1}},$$

$$I_3 = \sqrt{-x^2+x+1} - \int \frac{dx}{2\sqrt{-x^2+x+1}} \tag{18}$$

Let

$$I_4 = \int \frac{dx}{2\sqrt{-x^2+x+1}} \tag{19}$$

$$I_4 = \int \frac{dx}{\sqrt{-4x^2+4x+4}} = \int \frac{dx}{\sqrt{5-(1-2x)^2}} = \int \frac{\frac{1}{\sqrt{5}}\,dx}{\sqrt{1-\left(\frac{1-2x}{\sqrt{5}}\right)^2}} = -\frac{1}{2}\int \frac{\frac{-2}{\sqrt{5}}\,dx}{\sqrt{1-\left(\frac{1-2x}{\sqrt{5}}\right)^2}}$$

$$I_4 = -\frac{1}{2}arcsin\left(\frac{1-2x}{\sqrt{5}}\right) \tag{20}$$

From (18) and (19), $I_3 = \sqrt{-x^2+x+1} - I_4$; substituting in this the result given by (20),

$$I_3 = \sqrt{-x^2+x+1} + \frac{1}{2}arcsin\left(\frac{1-2x}{\sqrt{5}}\right) \tag{21}$$

From (16), $I_1 = x\,arcsin\left(\frac{1-2x}{\sqrt{5}}\right) - I_3$, and substituting therein the result given by (21),

$$I_1 = \left(x-\frac{1}{2}\right)arcsin\left(\frac{1-2x}{\sqrt{5}}\right) - \sqrt{-x^2+x+1} \tag{22}$$

And by a procedure completely analogously, we obtain

$$I_2 = \left(x+\frac{1}{2}\right)arcsin\left(\frac{1+2x}{\sqrt{5}}\right) + \sqrt{-x^2-x+1} \tag{23}$$

Substituting in (13) this results, we obtain the exact value of the requested probability:

$$p = \frac{1}{\pi}\left[I_1 + I_2\right]_{-\frac{1}{2}}^{\frac{1}{2}} + \frac{2}{\pi}arctan(2) - 1 = \frac{1}{\pi}\left[2arctan\left(\frac{1}{3}\right) + \frac{\pi}{2} + 1 - \sqrt{5}\right] + \frac{2}{\pi}arctan(2) - 1$$

$$\boxed{p = \frac{2}{\pi}\left[arctan\left(\frac{1}{3}\right) + arctan(2)\right] - \frac{\sqrt{5}-1}{\pi} - \frac{1}{2}} \quad (24)$$

The expression (24) can be simplified considering the definition of the *golden ratio* [1] and the following identity regarding tangent arcs (by the general shape established in [2] for the decomposition of pi / 4 in two arctan):

$$\boxed{arctan(2) = arctan\left(\frac{1}{3}\right) + \frac{\pi}{4}} \quad (25)$$

This identity can be proven easily by the formula of the tangent of a sum or through algebra of complex numbers, expressing the product of two complex numbers (suitably chosen) in two representation forms, binary form and polar form, as shown below.

Product in **binary form** and its corresponding representation in **polar form**:

$$(3+i)(1+i) = (2+4i) \Leftrightarrow \sqrt{10}_{arctan\left(\frac{1}{3}\right)} \sqrt{2}_{\frac{\pi}{4}} = \sqrt{20}_{arctan(2)}$$

After performed the product in polar form, the identity (25) is derived by identifying the arguments on both sides of the last equality:

$$\boxed{\sqrt{20}_{arctan\left(\frac{1}{3}\right)+\frac{\pi}{4}} = \sqrt{20}_{arctan(2)}} \quad (26)$$

Finally, the result (24) can be expressed in the following elegant form that involves two of the most remarkable numbers: the **number Pi** and the **Golden Ratio** Φ,

$$\boxed{\Phi = \frac{1+\sqrt{5}}{2}} \quad (27)$$

As the number π, it is surprising the ubiquity of this number, that emerge in the most diverse sceneries [1].

$$\boxed{p = \frac{2}{\pi}\left(2\arctan\frac{1}{3} - \frac{1}{\Phi}\right)} \quad (28)$$

The approximate value of *p* in ten thousandths is, *p* ≈ 0.0162.

## References.